\documentclass[12pt]{article}
\def\N{{\mathbf{N}}}
\def\Tr{{\mathrm{Tr}}}
\def\res{{\mathrm{res}}}
\def\bC{{\mathbf{\overline{C}}}}
\def\Ima{{\mathrm{Im}\,}}
\def\Rea{{\mathrm{Re}\,}}
\def\inter{{\mathrm{int}}}
\begin{document}
\title{Herbert Stahl's proof of the BMV conjecture}
\author{Alexandre Eremenko\footnote{Supported by NSF.}}
\maketitle

\begin{abstract} The paper contains a simplified version of
Stahl's proof of a conjecture of Bessis, Moussa and Villani
on the trace of matrices $A+tB$ with Hermitian $A$ and $B$.

MSC 2010 Class: 47A55, 30F10. Keywords: Hermitian matrices, perturbation
theory, trace, Riemann surfaces.
\end{abstract}

This paper presents a simplified
version of the proof of Herbert Stahl's theorem
on the BMV conjecture \cite{S}. The proof preserves
all main ideas of Stahl; the simplification consists in
technical details.
\vspace{.1in}

\noindent
{\bf Theorem.} {\em Let $A$ and $B$ be
two $n\times n$ Hermitian matrices. 
Then the function
\begin{equation}
\label{f}
f(t)=\Tr\;e^{A-tB}
\end{equation}
has a representation
\begin{equation}\label{mu}
f(t)=\int_{b_1}^{b_n} e^{-st}d\mu(s),
\end{equation}
where $\mu$ is a non-negative measure, $b_1$ and $b_n$
are the smallest and the largest eigenvalues of $B$.}
\vspace{.1in}

If $B$ is positive semi-definite, it follows that
$(-1)^nf^{(n)}\geq 0$. Such functions are called absolutely monotone.
The result was conjectured in \cite{BMV}.
Two equivalent statements for positive semi-definite matrices $B$ are
that
the polynomial $t\mapsto\Tr(A+Bt)^p,\; p\in\N,$
has all non-negative coefficients,
and that the function $t\mapsto\Tr(A+tB)^{-p},\; p\geq 0$ is absolutely
monotone,
\cite{LS}. Before the work of Stahl,
Theorem 1 was known for $2\times 2$ matrices. 
The proof of Stahl, which is explained in these notes,
is completely elementary: all needed tools were available in
the middle of XIX century.

Without loss of generality, one can assume that $B$ is a diagonal matrix
with eigenvalues $b_n>b_{n-1}>\ldots>b_1>0$. This is achieved by simultaneous
conjugacy of $A$ and $B$, adding a scalar to $B$,
and approximating the resulting $B$ with a matrix whose
eigenvalues are distinct.

Now eigenvalues $\lambda$ of $A-tB$ are determined from the equation
$$\det(\lambda I-A+tB)=0.$$
This determinant is a polynomial in two variables $t,\lambda$. We take $t$
out of the determinant, and denote $y=\lambda/t,\; x=1/t$, then
we obtain a polynomial equation of the form 
$$0=\det(yI+B-xA)=\prod_{j=1}^n(y+b_j-xa_{j,j})+O(x^2),$$
where $O(x^2)$ is a polynomial divisible by $x^2$.

This implies that there are $n$ holomorphic branches of the multivalued
implicit function $\lambda(t)$ in a neighborhood of infinity, which
satisfy
\begin{equation}
\label{branches}
\lambda_j(t)=-b_jt+a_{j,j}+O(1/t),\quad t\to\infty,
\end{equation}
and all $\lambda_j$ are real on the real line. Moreover, each of these branches
has an analytic continuation in a region containing the real line,
according to Rellich's theorem \cite[Thm XII.3]{RS}.
The algebraic function $\lambda(t)$ is defined on a Riemann
surface $S$ with $n$ sheets spread over the Riemann sphere.
This Riemann surface is not necessarily connected. It has $n$
unramified sheets over a region that contains
the real line and a neighborhood of infinity.
\vspace{.1in}

\noindent
{\bf Special case.} Suppose that $A$ is also diagonal, then
the $O(1/t)$ terms in (\ref{branches}) can be omitted,
and we obtain
\begin{equation}\label{special}
f(t)=\sum_{j=1}^n e^{a_{j,j}}e^{-b_jt}=\int_0^\infty e^{-st}
\sum_{j=1}^ne^{a_{j,j}}
\delta_{b_j}(s)ds.
\end{equation}
Thus $\mu$ is a discrete measure with positive atoms at the eigenvalues
of $B$. 
\vspace{.1in}

In the general case, the discrete component of $\mu$ is the same,
and the continuous component is a positive function on $(b_1,b_n)$.

Stahl figured out the following explicit
expression for the density.
\vspace{.1in}

\noindent
{\bf Proposition 1.} {\em The measure
\begin{equation}\label{density}
d\mu(s)=\left(\sum_{j=1}^n e^{a_{j,j}}\delta_{b_j}(s)+w(s)\right)ds
\end{equation}
where
\begin{equation}\label{w}
w(s)=-\sum_{j:b_j<s}\res_\infty e^{\lambda_j(\zeta)+s\zeta}=
\frac{1}{2\pi i}\sum_{j:b_j<s}\int_Ce^{\lambda_j(\zeta)+s\zeta}d\zeta,
\end{equation}
satisfies (\ref{f}) and (\ref{mu}). Here $C$ is any circle centered
at the origin, of sufficiently large radius, described counterclockwise.}
\vspace{.1in}

First we give an heuristic argument which could be used to guess this formula.
Inversion formula for the Laplace transform gives the density in the form
$$\frac{1}{2\pi i}\int_Lf(\zeta)e^{s\zeta}d\zeta,$$
where $L$ is a vertical line sufficiently far to the right.
For $|\zeta|$ large enough, the expression under the integral equals
$$e^{s\zeta}f(\zeta)=\sum_{j=1}^ne^{\lambda_j(\zeta)+s\zeta}.$$
As $\lambda_j(\zeta)=-b_j\zeta+\ldots,$ the summands for which $b_j>s$
are exponentially decreasing in the right half-lane, therefore, for these
summands the line $L$ can be shifted to the right, and all these summands 
vanish. The rest of the summands exponentially decrease to the left,
and for them, the contour can be bent to the left to obtain a circle $C$.

Of course one can give a rigorous justification of these arguments,
but once the formula is guessed, it is easy to verify it directly,
and we reproduce Stahl's argument.
\vspace{.1in}

\noindent
{\bf Lemma 1.} {\em For every $s$, we have
$$\sum_{j=1}^n\int_Ce^{\lambda_j(\zeta)+s\zeta}d\zeta=0.$$}

Indeed, this is an integral of an entire function over a closed contour.
\vspace{.1in}

It follows that the density $w$ defined by (\ref{w}) is zero for $s>b_n$,
and it is evidently zero for $s<b_1$.
\vspace{.1in}

{\em Proof of Proposition 1.} We compute the Laplace transform of
the density $w$ defined by (\ref{w}).
$$\int_0^\infty e^{-st}w(s)ds=\sum_{k=1}^{n-1}\int_{b_k}^{b_{k+1}}e^{-ts}w(s)ds=:\sum_{k=1}^{n-1}
I_k(t).$$
We fix $t>0$ and deform the contour $C$ in (\ref{w}) so that the positive ray
is outside $C$. This is possible to do because all $\lambda_j$
are holomorphic in a region containing the real line
and $C$. Thus $t$ is outside of the deformed contour $C^\prime$.
According to (\ref{w}), we have 
$$I_k(t)=\int_{b_k}^{b_{k+1}}\sum_{j=1}^k\frac{1}{2\pi i}\int_{C^\prime}
e^{\lambda_j(\zeta)+s(\zeta-t)}d\zeta ds.$$
By changing the order of integration and the order of summation,
we obtain
\begin{eqnarray*}
\sum_{k=1}^{n-1}I_k(t)&=&\sum_{j=1}^{n-1}\frac{1}{2\pi i}\int_{C^\prime}
e^{\lambda_j(\zeta)}\int_{b_j}^{b_n}e^{s(\zeta-t)}ds d\zeta\\
&=&\sum_{j=1}^{n-1}\frac{1}{2\pi i}\int_{C^\prime}e^{\lambda_j(\zeta)}\left(e^{b_n(\zeta-t)}-
e^{b_j(\zeta-t)}\right)\frac{d\zeta}{\zeta-t}.
\end{eqnarray*}
The last expression is transformed using Cauchy's formula and the fact
that $t$ is outside $C^\prime$. We have
$$\sum_{j=1}^n\int_{C^\prime}e^{\lambda_j(\zeta)}e^{b_n(\zeta-t)}
\frac{d\zeta}{\zeta-t}=0,$$
similarly to Lemma 1,
so 
$$\sum_{k=1}^{n-1}I_k(t)=
-\sum_{j=1}^n\frac{1}{2\pi i}\int_{C^\prime}e^{\lambda_j(\zeta)+b_j(\zeta-t)}\frac{d\zeta}{\zeta-t}.$$
Using (\ref{branches}), we write
$\lambda_j(\zeta)=-b_j\zeta+a_{j,j}+r_j(\zeta)$, where $r_j(\infty)=0$,
and apply Cauchy's formula again.
We obtain for every $j$:
\begin{eqnarray*}
&-&\frac{1}{2\pi i}\int_{C^\prime} e^{\lambda_j(\zeta)+b_j(\zeta-t)}\frac{d\zeta}{\zeta-t}
=
-\frac{e^{-b_jt+a_{j,j}}}{2\pi i}
\int_{C^\prime}e^{r_j(\zeta)}\frac{d\zeta}{\zeta-t}\\ \\
&=&
e^{-b_jt+a_{j,j}}\left(e^{r_j(t)}-1\right)=e^{\lambda_j(t)}-e^{-b_jt+a_{j,j}}.
\end{eqnarray*}
Adding these expressions for $j=1\ldots n$ and comparing with (\ref{density})
and the second equation in (\ref{special}),
we obtain Proposition 1.
\vspace{.1in}

It remains to prove that (\ref{w}) is non-negative for every $s$. 
Let us fix $s$ and $k$ so that $b_k<s<b_{k+1}$.
The idea of Stahl, is to replace the contour of integration in (\ref{w})
by an ingeniously chosen homologous contour,
on which the integral is non-negative
simply because the integrand is non-negative.

We recall that $S$ is a (possibly disconnected) Riemann surface
spread over the $\zeta$-sphere.
We denote a generic point of $S$ by $p$, and let $\pi:S\to\bC$
be the projection to $\zeta$-plane. Then $\lambda$ is a meromorphic function
on $S$  whose all poles are simple and lay over $\zeta=\infty$.

Asymptotic expressions (\ref{branches}) imply that there exists $R>0$
such that for all $j\leq k$ the functions
\begin{equation}\label{bra}
\lambda_j(\zeta)+s\zeta=(s-b_j)\zeta+\ldots
\end{equation}
are holomorphic for $|\zeta|>R$, real on the real line and 
have strictly positive derivatives for $\zeta>R$ and $\zeta<-R$,
while for $j>k$ they have strictly negative derivatives.
By increasing $R$, if necessary, we achieve
that for $|\zeta|>R/2$ and $j\leq k$, we have that
$\Ima (\lambda_j(\zeta)+s\zeta)$
has the same sign as $\Ima\zeta$.
And for $|\zeta|>R/2$ and $j>k$,
$\Ima (\lambda_j(\zeta)+s\zeta)$ has the opposite sign from $\Ima\zeta$. 

The surface $S$ has an anti-conformal involution, induced by complex
conjugation. The set of fixed points of this involution consists of
$n$ curves, $\pi$-preimages of the real line.
These curves break $S$ into two halves
$S^+$ and $S^-$ which are mapped onto each other by the involution.
Projections of these halves are the upper and lower half-planes.

We set $C=\{ \zeta:|\zeta|=R\}$ in (\ref{w}), where $R$ was just chosen. 

Consider the open sets
$$D^+:=\{ p\in S: |\pi(p)|<R,\; \Ima\pi(p)>0,\;\Ima(\lambda(p)+s\pi(p))>0\},$$
$$D^-:=\{ p\in S: |\pi(p)|<R,\; \Ima\pi(p)<0,\;\Ima(\lambda(p)+s\pi(p))<0\},$$
and 
$$D=\inter\left(\overline{D^+}\cup\overline{D^-}\right).$$
The set $\{ p\in S: |\pi(p)|=R\}$ consists of $n$
disjoint circles
$C_j\subset S$ which we label according to the branches
of $\lambda_j$ in (\ref{bra}), so that $\lambda=\lambda_j$ on $C_j$.
According to the paragraph after (\ref{bra}), the circles 
$C_j$ with $j\leq k$ belong to $\partial D$ while the $C_j$ with $j>k$ are
disjoint from $\overline{D}$. 

Let $D_1$ be a component of $D$ whose boundary contains some circles
$C_j$.\footnote{One can prove using the maximum principle that
every component of $D$ has some $C_j$ on the boundary, but
we are not using this fact.}
We are going to prove that 
\begin{equation}\label{pos}
\sum_{j:C_j\subset\partial D_1}\int_{C_j}e^{\lambda(p)+s\pi(p)}d\pi(p)>0,
\end{equation}
where the circles are oriented counterclockwise, which agrees
with their orientation as part of $\partial D$.
Adding these relations over all components of $D$
will prove the theorem.
Indeed, each circle $C_j$ with $j\leq k$
belongs to the boundary of exactly one component of $D$,
and circles $C_j$ with $j>k$ do not belong to the boundary of $D$.

Each component $D_1$ of $D$ is a Riemann surface of finite type, whose boundary
consists of several curves parametrized by circles. 
This parametrization is piecewise smooth, but may be
neither smooth nor injective.
We call these curves the boundary curves of $D_1$.
Our choice of $R$ guarantees that
the part of the boundary of $D_1$
that projects in $C$ is exactly the
chain on which the integration is performed in (\ref{pos}).
Consider the rest of the boundary $\partial D_1$
which projects into $|\zeta|<R$. 
\vspace{.1in}

\noindent
{\bf Lemma 2.} {\em No boundary curve of $D$ over $\{\zeta:|\zeta|<R\}$
can project into the open
upper or lower half-plane.}
\vspace{.1in}

Indeed suppose that $\gamma$ is a boundary curve whose projection does not
intersect the real axis. It is oriented in the standard way, so that $D$
is on the left. Suppose without loss of generality
that $\gamma$ projects to the upper half-plane.
Let $g(p)=\lambda(p)+s\pi(p)$. As $\Ima g>0$ in $D^+$,
and $\Ima g=0$ on $\gamma$, we conclude that the normal
derivative of $\Ima g$ has constant sign on $\gamma$. Then by the
Cauchy-Riemann
equations, the tangential derivative of $\Rea g$ along $\gamma$
is of constant sign, which is impossible because $\gamma$ is a closed curve,
and $\Rea g$ is single valued on $\gamma$.
\vspace{.1in}

Thus every boundary curve of $D_1$ intersects the real line.
Let $\gamma$ be a boundary curve of $D_1$ which projects into $\{\zeta:
|\zeta|<R\}$.
By Lemma 2, $\gamma$ is mapped into itself by the involution, so it consists
of two symmetric pieces: one piece $\gamma^+$ projects in the upper
half-plane, another $\gamma^-$ to the lower half-plane.
At all endpoints $p$ of
$\gamma^+$ or $\gamma^-$ we have $\Im\pi(p)=0$.  We have 
$$e^g=e^{\Rea g+i\Ima g}=e^{\Rea g}\left(\cos(\Ima g)+i\sin(\Ima g)\right).$$
Since $\Ima g=0$ on $\gamma$, and $\Rea g$ is increasing, we conclude that
$\phi(t):=e^{g(\gamma(t))}$ is real and increasing function
of the natural parameter $t$
on $\gamma^+$.
Thus
\begin{eqnarray*}
\frac{1}{2\pi i}\int_\gamma e^{g(p)}d\pi(p)&=&
\frac{1}{2\pi i}\left(\int_{\gamma^+}+\int_{\gamma^-}\right)
\phi(t)(d\xi(t)+id\eta(t))\\
&=&\frac{1}{\pi}\int_{\gamma^+}\phi(t)d\eta(t)
=
-\frac{1}{\pi}
\int_{\gamma^+}\eta(t)d\phi(t)<0,
\end{eqnarray*}
where we integrated by parts using $\eta(t)=0$ on the endpoints of
$\gamma^+$.

As the integral of the holomorphic $1$-form over the boundary equals zero,
$$\int_{\partial D_1}e^gd\pi=0$$
by Cauchy's theorem, the contribution to the integral from the
part of $\partial D_1$ which projects to $\{\zeta:|\zeta|<R\}$
is the negative of the contribution
of the part of $\partial D_1$ over $C$.
This completes the proof of (\ref{pos}) and of Theorem 1.

{\em Purdue University

eremenko@math.purdue.edu}

\end{document}